\documentclass[times,sort&compress,3p]{elsarticle}

\usepackage{amsmath,amsfonts,amssymb,amsthm,booktabs,color,epsfig,graphicx,hyperref,url}
\usepackage{tabularx}
\usepackage{amsthm}
\usepackage{here}
\usepackage{caption}
\usepackage{subcaption}
\usepackage{ascmac}
\usepackage{color}
\usepackage[labelfont=bf]{caption}
\usepackage{latexsym}
\usepackage{comment}

\newtheorem{theorem}{Theorem} 
\newtheorem{definition}{Definition}
\newtheorem{corollary}{Corollary}
\newtheorem{algorithm}{Algorithm} 
\newtheorem{example}{Example} 

\newtheorem{lemma}[theorem]{Lemma} 
\newdefinition{rmk}{Remark} 
\newproof{pf}{Proof}
\newproof{pot}{Proof of Theorem \ref{thm2}}

\begin{document}

\begin{frontmatter}



\title{Algorithm for the product of Jack polynomials and \\its application to the sphericity test}


\author[a]{Koki Shimizu}
\author[a]{Hiroki Hashiguchi}

\cortext[mycorrespondingauthor]{Corresponding author. Email address: \url{1420702@ed.tus.ac.jp}~(K. Shimizu).}

\affiliation[a]{organization={Tokyo University of Science},
            addressline={1-3 Kagurazaka}, 
            city={Shinjuku-ku},
            postcode={162-8601}, 
            state={Tokyo},
            country={Japan}}

\begin{abstract}
In this study, we derive the density and distribution function of a ratio of the largest and smallest eigenvalues of a singular beta-Wishart matrix for the sphericity test. 
These functions can be expressed in terms of the product of Jack polynomials. 
We propose an algorithm that expands the product of Jack polynomials by a linear combination of Jack polynomials. 
Numerical computation for the derived distributions is performed using the algorithm. 

\end{abstract}

\end{frontmatter}


\section{Introduction}

\noindent
The zonal polynomial, which is a real symmetric homogeneous polynomial, appears in the distribution theory of eigenvalues in multivariate analysis. 
The properties and essential integral formulas of zonal polynomials are given by Constantine~\cite{C1963} and James~\cite{James1964}.
Some extensions of zonal polynomials are useful in deriving distributions related to the eigenvalues of the central or non-central Wishart matrices. 
The product of two zonal polynomials are represented by a linear combination of zonal polynomials in Constantine~\cite{C1966} and Hayakawa~\cite{Hayakawa1967}.
 As regards multivariate analysis of variance (MANOVA), Khatri and Pillai~\cite{Khatri1968} and Pillai and Sugiyama~\cite{Pillai1969} proposed the non-null density function and distribution function for Roy's test statistic, respectively. 
 Sugiyama~\cite{Sugiyama1970} derived the density of the ratio of the largest and smallest eigenvalues of a nonsingular Wishart matrix for the sphericity test. 
 Ratnarajah et al.~\cite{Ratnarajah2004} extended this result to the complex case and applied the density to communication systems. 
 Their study essentially entails the use of Schur polynomials (called complex zonal polynomials) instead of zonal polynomials. 
 The aforementioned functions can be represented by the product of zonal or Schur polynomials. 
  However, in  terms of numerical calculations, this representation is not useful because algorithms for the product of zonal and Schur polynomials are not known. 
 
 Another extension of zonal polynomials is the Jack polynomials. 
The Jack polynomial, introduced by Jack~\cite{Jack1970}, is a symmetric homogeneous polynomial of which the zonal and Schur polynomials are special cases. D\'iaz-Garc\'ia and Guti\'errez-J\'aimez~\cite{Garcia2011} derived the density of a nonsingular beta-Wishart matrix, which covers the classical matrix of real, complex, and quaternion cases. 
D\'iaz-Garc\'ia~\cite{Garcia2014} derived the exact distribution of the largest and smallest eigenvalues of the nonsingular beta-Wishart matrix in terms of the Jack polynomials. 
Recently, Shimizu and Hashiguchi~\cite{Shimizu2021a} extended the distribution theory of eigenvalues of Wishart matrices from nonsingular to singular. 
The exact distributions of the largest eigenvalue of the singular beta $F$-matrix and elliptical Wishart matrix are given by Shimizu and Hashiguchi~\cite{Shimizu2021b} and Shinozaki et al.~\cite{Shinozaki2021}, respectively.

In this paper, we propose an algorithm to expand the product of Jack polynomials by using a linear combination of those polynomials. 
Furthermore, we derive the exact distribution of the ratio of the extreme eigenvalues of a singular beta-Wishart matrix for the sphericity test when the sample size is less than the dimension.
In Section~\ref{sec:02}, we define the Jack polynomials in terms of elementary symmetric functions and propose an algorithm for computing the products of Jack polynomials. 
The density and distribution functions for the ratio of extreme eigenvalues are obtained in a form that includes a linear combination of Jack polynomials in Section~\ref{sec:03}.
The derivation of their distribution is fundamentally the same as that of Sugiyama~\cite{Sugiyama1970}. 
Numerical computations of the distribution of the ratio of the extreme eigenvalues are presented in Section~\ref{sec:04}.

\section{Algorithm for product of Jack polynomials}\label{sec:02}
Let $\mathbb{F}_\beta$ denote a real finite-dimensional division algebra such that
$\mathbb{F}_1 = \mathbb{R}$, $\mathbb{F}_2 =\mathbb{C}$, and $\mathbb{F}_4 = \mathbb{H}$ for $\beta=1, 2, 4$, where
$\mathbb{R}$ and $\mathbb{C}$ are the fields of real and complex numbers, respectively, and $\mathbb{H}$ is the quaternion division algebra over $\mathbb{R}$.
Let $\mathbb{F}_\beta^{m \times n}$ be denoted by the set of all $m \times n$ matrices over $\mathbb{F}_\beta$, where $m \ge n$.
The conjugate transpose of $X \in \mathbb{F}_\beta^{m \times m}$ is written as $X^\ast=\overline{X^\top}$, and we say that $X$ is Hermitian if $X^\ast = X$.  The set of all Hermitian matrices is denoted by 
$S^{\beta}(m)=\{ X \in \mathbb{F}^{m\times m}_\beta \mid X^\ast=X\}$. 
The eigenvalues of the Hermitian matrix are real.  If the eigenvalues of $X \in S^\beta(m)$ are all positive, 
it can be considered positive definite and $X > 0$.
We can represent the Stiefel manifold over $\mathbb{F}_\beta$ as
$V^{\beta}_{n,m} = \{ H_1\in \mathbb{F}_\beta^{m\times n}\mid H_1^{\ast}H_1=I_n\}$.
For $c \in \mathbb{F}_\beta$, the multivariate gamma function of parameter $\beta>0$ is defined by
\begin{align}\label{def-gamma}
\nonumber
\Gamma_m^{\beta}(c) &= \int_{X > 0} |X|^{c - (m-1)\beta/2 } \mathrm{etr}(-X) (d X)
\\ 
 &= \pi^{\frac{m(m-1)\beta}{4}}\prod_{i=1}^{m}\Gamma\bigg\{c-\frac{(i-1)\beta}{2}\bigg\},
\end{align}
where $\Re(c)>\frac{(m-1)\beta}{2}$, $|X|$ are the determinant of a matrix $X$ and $\mathrm{etr}(\cdot)=\mathrm{exp}(\mathrm{tr}(\cdot))$.
 For a positive integer $k$, let $\kappa=(\kappa_1,\kappa_2,\dots,\kappa_m)$ denote a partition of $k$ with $\kappa_1\geq\cdots+\kappa_m\geq 0$ and $\kappa_1+\cdots +\kappa_m=k$. 
 The set of all partitions with lengths less than $m$ is denoted by $P^k_{m}=\{\ \kappa=(\kappa_1,\dots,\kappa_m)\mid \kappa_1+\dots+\kappa_m=k, \kappa_1\geq \kappa_2\cdots \geq \kappa_m \geq 0 \}$. 
The generalized Pochhammer symbol of parameter $\beta>0$ is defined by
\begin{align*}
(a)_\kappa^{\beta}=\prod_{i=1}^{m}\biggl(a-\frac{i-1}{2}\beta\biggl)_{\kappa_i}.
\end{align*}
The elementary symmetric functions in eigenvalues $x_1,  \dots, x_m$ of $X$ are expressed as $e_0=1$ and 
\begin{align*}
e_1=x_1+\cdots +x_m, ~~e_2=x_1x_2+\cdots+x_{m-1}x_m, ~~e_m=x_1x_2\cdots x_m.
\end{align*}
For $\kappa\in P^k_{m}$, we define the polynomials $\mathcal{E}_\kappa(X)$ as
\begin{align*}
\mathcal{E}_\kappa(X)=e_1^{\kappa_1-\kappa_2}e_2^{\kappa_2-\kappa_3}\cdots e_{m-1}^{\kappa_{m-1}-\kappa_{m}}e_m^{\kappa_m}, 
\end{align*}
where the degree of $\mathcal{E}_\kappa(X)$ is $(\kappa_1-\kappa_2)+2(\kappa_2-\kappa_3)+\cdots+m\kappa_m=m$.
The aforedescribed definition was also reported by Takemura~\cite{Takemura1984} and Jiu and Koutshan~\cite{Jiu2020}. 
Hashiguchi et al.~\cite{Hashiguchi2000} defined zonal polynomials, which is a special case of Jack polynomials in terms of elementary symmetric functions. 
Similar to the description by Hashiguchi et al.~\cite{Hashiguchi2000}, the definition of Jack polynomials in terms of elementary symmetric functions can be stated as follows.
\begin{definition}\label{def1} 
 For $\kappa \in P^k_m$ and $X\in S^\beta(m)$, there exists a unique polynomial $C^\beta_\kappa(X)$ that satisfies the following three conditions: 
\begin{itemize}
\item[1] $C^\beta_\kappa(X)=\sum_{\mu\preceq \kappa} q^\beta[\kappa, \mu]\mathcal{E}_\mu(X)$ and $q^\beta[\kappa, \kappa] \neq 0$,
where $q^\beta[\kappa, \mu]$ is a constant, and the summation is over all partitions $\mu\preceq \kappa$; that is, $\mu$ is equal to or less than $\kappa$ in the lexicographical ordering. 

 \item[2] $D^\beta_{m}C^\beta_\kappa(X)=d^\beta(\kappa) C^\beta_\kappa(X)$,
 where $D^\beta_{m}$ is the Laplace--Beltrami operator
 $ D^\beta_{m}=
    \sum_{i=1}^{m}x_{i}^{2}\frac{\partial^{2}}{\partial x_{i}^{2}}
    +\beta \sum_{i\leq i \neq j \leq m}\frac{x_{i}^{2}}{x_{i}-x_{j}}\frac{\partial}{\partial x_{i}}$, $d^\beta(\kappa)$ is the corresponding eigenvalue. 
    $d^\beta(\kappa)=\sum_{i=1}^{m}\lambda_i(\beta\lambda_i+m-i-\beta)$.
    
 \item[3] $[\mathrm{tr}(X)]^k=\sum_{\kappa\in P^k_m}C^\beta_\kappa (X)$.
\end{itemize}
 We hypothesize that the aforementioned polynomials $C^\beta_\kappa(X)$ are
 Jack polynomials. 
\end{definition}
An ordinal definition of Jack polynomials in terms of monomial symmetric functions is presented in Dumitriu et al.~\cite{Dumitriu2007}. 
Definition~\ref{def1} is useful in developing an algorithm for the product of Jack polynomials.
Hashiguchi et al.~\cite{Hashiguchi2000} determined the recurrence relation of $q^\beta[\kappa, \mu]$ for $\beta=1$.
It is easy to improve for any $\beta>0$ according to the same way of Hashiguchi et al.~\cite{Hashiguchi2000}.
\begin{corollary}\label{corollary1}
The coefficients $q^\beta[\kappa, \mu]$ satisfy the recurrence
 \begin{align}
 \label{coefficient-Jack}
 q^\beta[\kappa, \mu]=\frac{1}{d^\beta(\kappa)-d^\beta(\mu)}\sum_{\mu\preceq\nu\preceq\kappa}b^\beta[\nu, \mu]q^\beta[\kappa, \nu],
  \end{align}
 where $b^\beta[\nu, \mu]$ is determined by
 \begin{align*}
 D^\beta_{m}\mathcal{E}_\nu(X)=\sum_{\mu\preceq\nu}b^\beta[\nu, \mu]\mathcal{E}_\mu(X).
 \end{align*}
\end{corollary}
If $\kappa=\mu$, the coefficients $q^\beta[\kappa, \kappa]$ are represented by
\begin{align}
q^\beta[\kappa, \kappa]=\frac{(2/\beta)^k k!}{\prod_{(i, j)\in \kappa} h_\ast^\kappa(i, j)},
\label{qkk}
\end{align} 
where $h_\ast^\kappa(i, j)$ is the upper hook length at $(i, j)\in \kappa$. The study by Dumitriu et al.~\cite{Dumitriu2007} is also relevant in this regard. 

The key to the expansion of Jack polynomials $C^\beta_\kappa(X)$ in terms of $\mathcal{E}_\kappa(X)$ is the change in the variables of the Laplace--Beltrami operator from $x_1,\dots, x_m$ to $e_1,\dots, e_m$.
The recurrence relation \eqref{coefficient-Jack} for $m=2, 3, 4$ is provided in the Appendix.

The Jack polynomials for $\beta=1$ and $\beta=2$, are referred to as zonal polynomials and Schur polynomials, respectively. 
Herein, we represent $C^\beta_\kappa(X)$ for $\beta=1$ as $C_\kappa(X)$. 
For $\kappa\in P^k_m$, $\tau \in P^t_m$, the product of the Jack polynomials can be written as 
\begin{align}
\label{jackprod}
C^\beta_\kappa(X)\cdot C^\beta_\tau(X)=\sum_{\delta\in P^{k+t}_m}g^\delta_{\kappa, \tau}C^\beta_\delta(X), 
\end{align} 
where $g^\delta_{\kappa, \tau}$ is the linearization coefficient of $C^\beta_\delta(X)$.
If $\beta=1$ in (\ref{jackprod}), the coefficients  $g^\delta_{\kappa, \tau}$ for small $k$ and $t$ are calculated by Khatri and Pillai~\cite{Khatri1968} and Hayakawa~\cite{Hayakawa1967}. 
The following algorithm converts $\mathcal{E}_\kappa(X)$ into a linear combination of Jack polynomials $C_\kappa^\beta(X)$. 
This is a typical technique in computer algebra based on lexicographical ordering.
\begin{algorithm}{\label{algorithm1}} \rm
$$\text{Input: } \mathcal{E}_\nu(X) \qquad \text{Output: } \text{A linear combination of } \{C^\beta_{\kappa}(X)\},
$$
\begin{enumerate}
\item Set $r:=0$ and  $f := \mathcal{E}_\nu(X)$. 
\item If $f=0$, stop after returning $r$. 
\item Let $\kappa$ be the partition for the leading term of $f$ with respect to a linear combination of $\{\mathcal{E}_\mu (X)\}$.
We calculate all coefficients $q^\beta[\kappa, \mu]$ based on 
\eqref{coefficient-Jack} and \eqref{qkk} in order obtain get the identity~1 in Def.~\ref{def1} for $C^\beta_\kappa(X)$, as follows: 
$$
C^\beta_\kappa(X)=\sum_{\mu\preceq \kappa} q^\beta[\kappa, \mu]\mathcal{E}_\mu(X)
$$

\item Update $r$ and $f$ as
\begin{eqnarray*}
r:=r+\frac{\textrm{LC}(f)}{q^\beta[\kappa, \kappa]} \cdot C_{\kappa}^{\beta}(X),  
\quad 
f:= f -\frac{\textrm{LC}(f)}{q^\beta[\kappa, \kappa]} \cdot \left(
 \sum_{\mu\preceq \kappa} q^\beta[\kappa, \mu]\mathcal{E}_\mu(X)
\right). \\
 \end{eqnarray*}
 where $\textrm{LC}(f)$ is the leading coefficient of $f$, that is, the coefficient of $\mathcal{E}_\kappa(X)$ in $f$.
 Go back to 2.
\end{enumerate}
\end{algorithm}
\begin{example}
For partition $\kappa=(2,1, 1)$, the polynomials $\mathcal{E}_\kappa(X)$ can be expressed by
\begin{align*}
\mathcal{E}_{(2,1,1)}(X)=\frac{3}{16}C_{(2, 1, 1)}(X)+\frac{1}{2}C_{(1, 1, 1, 1)}(X).
\end{align*}
\end{example}
\begin{algorithm}\label{algorithm2} \rm
$$
\text{Input: } \mathcal{E}_{\mu}(X)\cdot \mathcal{E}_{\tau}(X) \qquad \text{Output: } \mathcal{E}_{\nu}(X)
$$
\begin{enumerate}
\item Let 
$$ \mathcal{E}_\mu(X) = \prod_{i=1}^m e_i^{\mu_{i} - \mu_{i+1}} \text{ and }. 
\mathcal{E}_\tau(X) = \prod_{i=1}^m e_i^{\tau_{i} - \tau_{i+1}},
$$
where $\mu_{m+1} = \tau_{m+1} = 0$ for convenience.
\item Set $\nu = (\mu_1 + \tau_1, \dots, \mu_m +\tau_m)$ and 
$$
 \mathcal{E}_\nu(X) = \prod_{i=1}^m e_i^{\nu_{i} - \nu_{i+1}},
$$
where $\nu_i = \mu_i  + \tau_i$ for $i=1, \dots, m+1$.
\item Return $ \mathcal{E}_\nu(X)$ as the product $\mathcal{E}_{\mu}(X)\cdot \mathcal{E}_{\tau}(X)$.
\end{enumerate}
\end{algorithm}

 By using Algorithms~\ref{algorithm1} and \ref{algorithm2}, the product of the Jack polynomials can be expressed as a linear combination of the Jack polynomials as follows:
\begin{algorithm}\label{algorithm3} \rm
$$\text{Input: } C^\beta_{\kappa}(X)\cdot C^\beta_{\tau}(X) \qquad \text{Output: }  \sum_{\delta\in P^{k+t}_m} g^\delta_{\kappa, \tau} C^\beta_{\delta}(X).
$$
\begin{enumerate}
\item Expand $C^\beta_\kappa(X)$ and $C^\beta_\tau(X)$ in terms of $\{ \mathcal{E}_\mu(X)\}$, respectively, based on identity~1 of Def.~\ref{def1}.
\begin{align}
 C^\beta_\kappa(X) &= \sum_{\mu\preceq \kappa} q^\beta[\kappa, \mu]\mathcal{E}_\mu(X)
 \label{prodCC-eq01}
 \\
C^\beta_\tau(X)  &=\sum_{\mu\preceq \tau} q^\beta[\tau, \mu]\mathcal{E}_\mu(X)
\label{prodCC-eq02}
\end{align}
\item Calculate the product of \eqref{prodCC-eq01} and \eqref{prodCC-eq02}, and apply Algorithm~\ref{algorithm2} to the right hand side of this product.
\item Obtain the above product as a linear combination of $\{\mathcal{E}_\mu \mid P_{m}^{k + t} \}$ and set it to $f$.
Apply Algorithm~\ref{algorithm1} to each term in $f$, and return 
$$
 \sum_{\delta\in P^{k+t}_m} g^\delta_{\kappa, \tau} C^\beta_{\delta}(X).
$$
\end{enumerate}
\end{algorithm}
\begin{example}
For partitions $\kappa=(2, 1), \tau=(2)$, $m=2$, the product of the Jack polynomials of $\beta=1$ can be expressed as
\begin{align*}  
C_{(2,1)}(X)\cdot C_{(2)}(X)
&=\frac{12}{5}\mathcal{E}_{(2, 1)}(X)\cdot \biggl(\mathcal{E}_{(2)}(X)-\frac{4}{3}\mathcal{E}_{(1, 1)}(X)\biggl)\\
&=\frac{12}{5}\mathcal{E}_{(4,1)}(X)-\frac{16}{5}\mathcal{E}_{(3, 2)}(X)\\
&=\frac{28}{75}C_{(3,2)}(X)+\frac{27}{50}C_{(4,1)}(X).
\end{align*}
\end{example}
Khatri and Pillai~\cite{Khatri1968} reported the product of zonal polynomials for some pairs of $\kappa$ and $\tau$ but encountered a few errors.
For example, we have $C_{(5)}(X) \cdot C_{(1)}(X) = C_{(6)}(X) + 5/27\; C_{(5, 1)}(X)$; the second coefficient was $1/54$ in  Khatri and Pillai~\cite{Khatri1968}.
\noindent

\section{Distribution of ratio of the largest and smallest eigenvalues for a singular beta Wishart matrix }\label{sec:03}
In this section, we discuss the ratio of the extreme eigenvalues of a singular beta-Wishart matrix for the sphericity test.
Let an $m\times n$ beta-Gaussian random matrix $X$ be distributed as $X\sim N_{m, n}^\beta(O, \Sigma \otimes \Theta)$, where $\Sigma$ and $\Theta$ are positive definite matrices. 
The density of $X$ is given as 
\begin{eqnarray*}
 \frac{1}{(2\pi \beta^{-1})^{mn\beta/2}|\Sigma|^{\beta n/2}|\Theta|^{\beta m/2}}\mathrm{exp}\biggl(-\frac{\beta}{2}\mathrm{tr}\Sigma^{-1}(X-M)\Theta^{-1}(X-M)^{\ast}\biggl).
\end{eqnarray*}
Let $X\sim N_{m, n}^\beta(O, \Sigma \otimes I_n)$. 
Then, the beta-Wishart matrix is defined as $W=XX^\ast$, and its distribution is denoted by $W_m^\beta(n, \Sigma)$. 
 If $n\geq m$, the random matrix $W$ is nonsingular; otherwise, it is singular.
 The spectral decomposition of the singular beta-Wishart matrix $W$ is $W=H_ 1 L _1H_1^\top$, where $L_1=\mathrm{diag}(\ell_1,\dots, \ell_n)$ and $\ell_1>\ell_2>\dots>\ell_n>0$.
 The density of $W$ for $\beta=1, 2$ was reported by Uhlig~\cite{U1994} and Ratnarajah et al. \cite{Ratnarajah2005}, respectively. 
Shimizu and Hashiguchi~\cite{Shimizu2021a} provided the density of $W$ as 
 \begin{align*}
  f(W)=\frac{\pi^{n(n-m)\beta/2}|\Sigma|^{-n\beta/2}}{(2\beta^{-1})^{mn\beta/2}\Gamma^{\beta}_n(n\beta/2)}|L_1|^{(n-m+1)\beta/2-1}\mathrm{etr}\biggl(-\frac{\beta}{2}\Sigma^{-1}W\biggl). 
 \end{align*}
We consider the sphericity test as  
\begin{align} \label{sphericity-test}
H_0:\Sigma=\sigma^2 I_m, \text{vs.}~H_1:\Sigma \neq \sigma^2 I_m,
\end{align}
where $\sigma^2$ is an unknown positive constant. 
For the sphericity test (\ref{sphericity-test}), the likelihood ratio test is usually performed for a fixed $m$ and a sufficiently large $n$. 
In contrast, if $m>n$, the likelihood ratio is not applicable because the sample covariance matrix is singular. 
Sugiyama~\cite{Sugiyama1970} proposed the ratio of extreme eigenvalues of the real nonsingular Wishart matrix as the statistic for (\ref{sphericity-test}) and obtained its density function. 
We extend these results to the singular case. 
The following lemma for $\beta=1$ was given by Sugiyama~\cite{Sugiyama1967}.
Shimizu and Hashiguchi~\cite{Shimizu2021a} generalized the result with parameter $\beta=1, 2, 4$ as follows:
  \begin{lemma} \label{lemma1}
  Let $X_1=\mathrm{diag}(1,x_2,\dots, x_n)$ and $X_2=\mathrm{diag}(x_2,\dots, x_n)$ with $x_2>\cdots >x_n>0$, then the following equation holds: 
  \begin{align*}
\int_{1>x_2>\cdots >x_n>0}|X_2|^{a-(n-1)\beta/2-1}C^{{\beta}}_\kappa(X_1)\prod_{i=2}^{n}(1-x_i)^\beta \prod_{i<j}(x_i-x_j)^\beta \prod_{i=2}^{n}dx_i\\
=(na+k)\{\Gamma^{{\beta}}_n(n\beta/2)/\pi^{n^2\beta/2+r_1}\}\frac{\Gamma^{{\beta}}_n(a,\kappa)\Gamma^{{\beta}}_n\{(n-1)\beta/2+1\}C^\beta_\kappa(I_n)}{\Gamma^{{\beta}}_n\{a+(n-1)\beta/2+1,\kappa\}},
\end{align*}
where $\mathrm{\Re}(a)>(n-1)\beta/2$, $\Gamma^{{\beta}}_n(\alpha,\kappa)=(\alpha)_\kappa \Gamma^\beta_n(\alpha)$, $\beta=1, 2, 4$, and 
\begin{align*}r_1&=
\begin{cases}
\quad 0,\quad  \quad~~~~\text{\ $\beta=1,$}\\
\quad -n\beta/2,~~~\text{$\beta=2,4$}.\\
\end{cases}
\end{align*}
\end{lemma}
\begin{theorem}\label{theorem1}
Let $W\sim W^\beta_m(n,\sigma^2I_m)$ with $m>n$. Then, the density of $x=1-\ell_n/\ell_1$ is given by
\begin{align}\label{maxmindensity}
\nonumber
f(x)&=C~\sum_{k=0}^{\infty}\sum_{\kappa\in P^k_{n-1}}\{\Gamma(mn\beta/2+k)/n^k k!\}\sum_{t=0}^{\infty}\{(n-1)(n\beta+2)/2+k+t\}/t!\\ 
&x^{(n-1)(n\beta+2)/2+k+t-1}\sum_{\tau\in P^{t}_{n-1}}\sum_{\delta\in P^{k+t}_{n-1}}\frac{g^\delta_{\kappa, 
\tau}\{(n-m-1)\beta/2+1\}^\beta_\tau (n\beta/2+1)^\beta_\delta C^\beta_\delta(I_{n-1})}{\{(n-1)\beta+2\}^\beta_\delta},
\end{align}
where $C=\frac{\Gamma^\beta_{n-1}(n\beta/2+1)\Gamma^\beta_{n-1}\{(n-2)\beta/2+1\}\pi^{r}}{\Gamma^\beta_{n-1}\{(n-1)\beta+2\}n^{mn\beta/2}\Gamma(n\beta/2)\Gamma^\beta_n(m\beta /2)}$ and $g^\delta_{ \kappa, \tau}$ is the coefficient of $C^\beta_\delta(Q)$, and
 \begin{align*}r&=
\begin{cases}
\quad n\beta/2, ~~~\text{$\beta=1,$}\\
\quad (n-1)\beta/2,~~~\text{$\beta=2,4$}.\\
\end{cases}
\end{align*}.
\begin{proof}
This proof is similar to that of Sugiyama~\cite{Sugiyama1970}.
From the joint density of the eigenvalues of a singular beta-Wishart matrix given by Shimizu and Hashiguchi~\cite{Shimizu2021a}, 
the null distribution for (\ref{sphericity-test}) is expressed as 
 \begin{align*}
 f(\ell_1,\dots,\ell_n)
 &=C_1~|L_1|^{(m-n+1)\beta/2-1}  \prod_{i<j}^{n}(\ell_i-\ell_j)^\beta \mathrm{etr}\biggl(-\frac{\beta}{2\sigma^2}\mathrm{tr}L_1\biggl),
\end{align*}
 where $C_1=\frac{(2\beta^{-1})^{-nm\beta/2}\pi^{n^2\beta/2+r_1}}{\sigma^{\beta mn}\Gamma^{{\beta}}_n(\frac{n\beta}{2})\Gamma^{{\beta}}_n(\frac{m\beta}{2})}$ and $r_1$ is given in Lemma~\ref{lemma1}.
Let $q_i=(\ell_1-\ell_i)/\ell_1$, for $i=2,\dots,n$.
Then, the joint densities of $\ell_1$ and $q_2,\dots, q_n$ is expressed as 
\begin{align*}
 f(\ell_1,q_2,\dots,q_n)=C_1~\mathrm{exp}\biggl(-\frac{\beta}{2\sigma^2}n\ell_1\biggl)|Q|^\beta \prod_{i>j}^{n}(q_i-q_j)^\beta|I_{n-1}-Q|^{(m-n+1)\beta/2-1} \sum_{k=0}^{\infty}\sum_{\kappa\in P^k_{n-1}}\ell_1^{mn\beta/2+k-1}\frac{C^\beta_\kappa(\frac{\beta}{2\sigma^2}Q)}{k!},
\end{align*}
where $Q=\mathrm{diag}(q_n,\dots,q_2)$ and $1>q_n>\cdots>q_2>0$.
Observe that
\begin{itemize}
\item[$1$.] $\mathrm{etr}\biggl(-\frac{\beta}{2\sigma^2}L_1\biggl)=\mathrm{exp}(-\frac{\beta}{2\sigma^2}n\ell_1)\sum_{k=0}^{\infty}\sum_{\kappa\in P^k_{n-1}}\frac{C^\beta_\kappa(\frac{\beta}{2\sigma^2}\ell_1Q)}{k!}$.
\item[$2$.] $\prod_{i<j}(\ell_i-\ell_j)^\beta=\ell_1^{n(n-1)\beta/2}|Q|^\beta \prod_{i>j}(q_i-q_j)^\beta$.
\item[$3$.] $(\mathrm{det}L_1)^{(m-n+1)\beta/2-1}=\ell_1^{ n(m-n+1)\beta/2-n}|I_{n-1}-Q|^{(m-n+1)\beta/2-1}$.
\end{itemize}
We use the fact that 
\begin{align*}
|I_{n-1}-Q|^{(m-n+1)\beta/2-1} C^\beta_\kappa(Q)&={_1F^{(\beta;n-1)}_0}\biggl(\frac{(n-m-1)\beta}{2}+1;Q\biggl)C^\beta_\kappa(Q)\\
&=\sum_{t=0}^{\infty}\sum_{\tau \in P^t_{n-1}}\frac{\{(n-m-1)\beta/2+1\}^\beta_ \tau C^\beta_\kappa(Q)C^\beta_ \tau(Q)}{t!}\\
&=\sum_{t=0}^{\infty}\sum_{ \tau \in P^{t}_{n-1}}\sum_{ \delta \in P^{k+t}_{n-1}}\frac{g^\delta_{\kappa, \tau}\{(n-m-1)\beta/2+1\}^\beta_ \tau C^\beta_\delta(Q)}{t!}.
\end{align*}
Translating $q_i$ to $s_i=q_i/q_n$, $i=2,\dots, n-1$, and using Lemma~\ref{lemma1}, 
we have 
\begin{align*}
&f(\ell_1,q_n)\\
&=C_1~\mathrm{exp}\biggl(-\frac{\beta}{2\sigma^2}n\ell_1\biggl)\sum_{k=0}^{\infty}\sum_{\kappa\in P^k_{n-1}}\frac{\ell_1^{mn\beta+k-1}(\beta/2\sigma^2)^k}{ k!}\sum_{t=0}^{\infty}\sum_{ \tau \in P^{t}_{n-1}}\sum_{ \delta \in P^{k+t}_{n-1}}\frac{g^\delta_{\kappa, \tau}\{(n-m-1)\beta/2+1\}^\beta_\tau}{t!}q_n^{n(n-1)\beta/2+k+t +n-2}\\
&\int_{1>s_{n-1}>\dots>s_2>0}|S|^\beta\prod_{i=2}^{n-1}(1-s_i)^\beta \prod_{i>j}(s_i-s_j)^\beta  \prod_{i=2}^{n-1} C^\beta_\delta(S) \prod _{i=2}^{n-1}ds_i\\
&=C_1~\mathrm{exp}\biggl(-\frac{\beta}{2\sigma^2}n\ell_1\biggl)\sum_{k=0}^{\infty}\sum_{\kappa\in P^k_{n-1}}\frac{\ell_1^{mn\beta+k-1}(\beta/2\sigma^2)^k}{ k!}\sum_{t=0}^{\infty}\sum_{ \tau \in P^{t}_{n-1}}\sum_{ \delta \in P^{k+t}_{n-1}}\frac{g^\delta_{\kappa, \tau}\{(n-m-1)\beta/2+1\}^\beta_\tau}{t!}q_n^{n(n-1)\beta/2+k+t+n-2}\\
&\{(n-1)(n\beta+2)/2+k+t\}\frac{\Gamma^\beta_{n-1}\{(n-1)\beta/2\} \Gamma^\beta_{n-1}(n\beta/2+1,\delta)\Gamma^\beta_{n-1}\{(n-2)\beta/2+1\}C^\beta_\delta(I_{n-1})}{\pi^{(n-1)^2\beta/2+r_2}\Gamma^\beta_{n-1}\{(n-1)\beta+2,\delta\}},
\end{align*}
where $S=\mathrm{diag}(s_{n-1},\dots,s_2)$, and $1>s_{n-1}>\cdots>s_2>0$ and 
\begin{align*}r_2&=
\begin{cases}
\quad 0,\quad  \quad~~~~\text{\ $\beta=1,$}\\
\quad -(n-1)\beta/2,~~~\text{$\beta=2,4$}.\\
\end{cases}
\end{align*}
From (\ref{def-gamma}), we have  
\begin{align*}
\frac{\Gamma^\beta_{n-1}\{(n-1)\beta/2\}}{\Gamma^\beta_n(n\beta/2)}=\frac{1}{\pi^{(n-1)\beta/2}\Gamma(n\beta/2)}.
\end{align*}
See also Eq. (3.5) in Kan and Koev~\cite{Kan2019}.
Finally, by integrating $f(\ell_1,q_n)$ with respect to $\ell_1$, we note the following identity: 
\begin{align*}
\int_{0}^{\infty}\ell_1^{mn\beta/2+k-1}\exp\biggl(-\frac{\beta}{2\sigma^2}n\ell_1 \biggl)d\ell_1=\frac{\Gamma(mn\beta/2+k)}{(n\beta/2\sigma^2)^{mn\beta/2+k}}, 
\end{align*}
and get (\ref{maxmindensity}).
\end{proof}
\end{theorem}
From Theorem~\ref{theorem1}, we can assume without loss of generality that  $\sigma^2=1$ for the null distribution of $\ell_n/\ell_1$. 
The product of the Jack polynomials in (\ref{maxmindensity}) can be calculated using Algorithm~\ref{algorithm2}. 
The generalization of the density of the ratio of extreme eigenvalues for nonsingular and singular beta-Wishart matrices for $\beta=1$ is given by 
\begin{align}\label{generalization-density}
\nonumber
f(x)&=C_2~\sum_{k=0}^{\infty}\sum_{\kappa\in P^k_{n_1-1}}\{\Gamma(n_1n_2/2+k)/n_1^k k!\}\sum_{t=0}^{\infty}\{(n_1-1)(n_1+2)/2+k+t\}/t!\\ 
&x^{(n_1-1)(n_1+2)/2+k+t-1}\sum_{\tau\in P^{t}_{n_1-1}}\sum_{\delta\in P^{k+t}_{n_1-1}}\frac{g^\delta_{\kappa, 
\tau}\{(n_1-n_2+1)/2\}_\tau \{(n_1+2)/2\}_\delta C_\delta(I_{n_1-1})}{ (n_1+1)_\delta},
\end{align}
where $n_1=\mathrm{min}(n, m)$, $n_2=\mathrm{max}(n, m)$ and $C_2=\frac{\pi^{n_1/2}\Gamma_{n_1-1}\{(n_1+2)/2\}\Gamma_{n_1-1}(n_1/2)}{\Gamma_{n_1-1}(n_1+1)n_1^{n_1n_2/2}\Gamma(n_1/2)\Gamma_{n_1}(n_2/2)}$. 
If $n\geq m$, the function (\ref{generalization-density}) coincides with the results obtained by Sugiyama~\cite{Sugiyama1970}. 
\begin{corollary}
Let $W\sim W^\beta_m(n,I_m)$, with $m>n$. 
The probability density function of $x=1-\ell_n/\ell_1$ is given as  
\begin{align}
\label{prob-ell1}
\nonumber
F(x)&=C~\sum_{k=0}^{\infty}\sum_{\kappa\in P^k_{n-1}}\{\Gamma(mn\beta/2+k)/n^k k!\}\sum_{t=0}^{\infty}x^{(n-1)(n\beta+2)/2+k+t}/t!\\
&\sum_{\tau \in P^{t}_{n-1}}\sum_{\delta\in P^{k+t}_{n-1}}\frac{g^\delta_{\kappa, \tau}\{(n-m-1)\beta/2+1\}^\beta_\tau (n\beta/2+1)^\beta_\delta C^\beta_\delta(I_{n-1})}{\{(n-1)\beta+2\}^\beta_\delta},
\end{align}
\begin{proof}
This proof is easily obtained from 
\begin{align*}
\int_{0}^{1}x^{(n-1)(n\beta+2)/2+k+t-1} dx=1/{\{(n-1)(n\beta+2)/2+k+t\}}x^{(n-1)(n\beta+2)/2+k+t}.
\end{align*}
\end{proof}
\end{corollary}
We also obtain the $h$-th moment of $1-\ell_n/\ell_1$ as 
\begin{align}
\label{moment}
\nonumber
E(x^h)&=C~\sum_{k=0}^{\infty}\sum_{\kappa\in P^k_{n-1}}\{\Gamma(mn\beta/2+k)/n^k k!\}\sum_{t=0}^{\infty}\biggl\{1-\frac{h}{(n-1)(n\beta+2)/2+k+t+h}\biggl\}/t!\\
&\sum_{\tau\in P^{t}_{n-1}}\sum_{\delta\in P^{k+t}_{n-1}}\frac{g^\delta_{\kappa, \tau}\{(n-m-1)\beta/2+1\}^\beta_\tau (n\beta/2+1)^\beta_\delta C^\beta_\delta(I_{n-1})}{\{(n-1)\beta+2\}^\beta_\delta}.
\end{align}

\section{Numerical experiments}\label{sec:04}
\begin{figure}[H]
\begin{center}
\includegraphics[width=7cm]{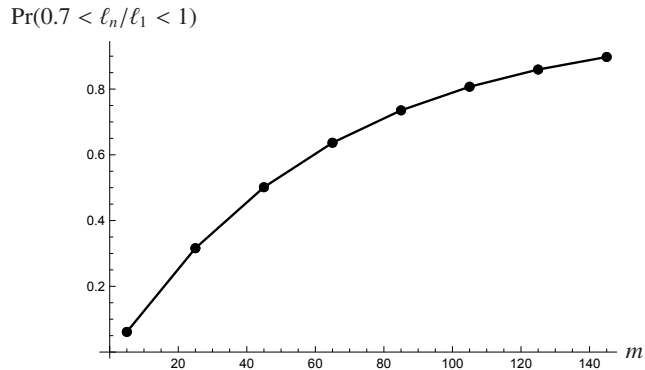}
\rlap{\raisebox{29.0ex}{\kern-23.em{\small $\Pr(0.7< \ell_n/\ell_1 < 1)$}}}%
\rlap{\raisebox{.15cm}{\kern0cm{\small $m$}}}
\caption{$n=2, \beta=1$}
 \label {fig1}
\end{center}
\end{figure}

This section presents the numerical computations performed for the derivation of the results.
If $p=\beta(m-n+1)/2-1$ is a positive integer, the truncated distribution up to the $K$ th degree of (\ref{prob-ell1}), which is a finite series for the summation of $t$, is represented by 
\begin{align*}
F_K(x)=&~C~\sum_{k=0}^{K}\sum_{\kappa\in P^k_{n-1}}\{\Gamma(mn\beta/2+k)/n^k k!\}\sum_{t=0}^{p(n-1)}x^{(n-1)(n\beta+2)/2+k+t-1}/t!\\
&\sum_{\tau^\ast}\sum_{\delta\in P^{k+t}_{n-1}}\frac{g^\delta_{\kappa, \tau}\{(n-m-1)\beta/2+1\}^\beta_\tau (n\beta/2+1)^\beta_\delta C^\beta_\delta(I_{n-1})}{\{(n-1)\beta+2\}^\beta_\delta},
\end{align*}
where $\sum_{\tau^\ast}$ is the sum of all partitions of $t$ with $\tau_1\leq p$ and $\tau=(\tau_1, \tau_2,\dots, \tau_{n-1})$.
The empirical distribution based on $10^6$ Monte Carlo simulations is denoted by $F_\mathrm{sim}$. 
We compute the percentage points of (\ref{prob-ell1}) as a finite series for $t$ in real~($\beta=1$) and complex $(\beta=2)$ cases. 
Table~\ref{table1} indicates the comparison $\alpha$ percentage points of $F_K(x)$ and $F_\mathrm{sim}$. 
We confirm that almost all percentage points achieve the desired accuracy.
Table~\ref{table1}~(a)--(b) indicates that the numerical computation of $F_K(x)$ in the case of a complex requires more series terms than the real cases. 
\begin{table}[H]
\caption{Percentile points of truncated distribution for (\ref{prob-ell1}).} \label{table1}
\begin{center}
\begin{tabular}{c}
    \begin{minipage}[c]{0.4\hsize}
      \begin{center}
       \captionsetup{labelformat=empty,labelsep=none}
         \subcaption{$\beta=1, m=10,n=3$}
         \label{table1-a}
{\begin{tabular}{@{}cccc@{}} \toprule
$\alpha$&${{F^{-1}_\mathrm{sim}}}(\alpha)$ &$F^{-1}_{25}(\alpha)$  \\ \toprule
0.01	 &~0.390& 0.389\\
0.05	 &~0.509& 0.509\\
0.50	 &~0.759& 0.759\\
0.90	 &~0.885 & 0.888\\ 
0.95	 &~0.910& 0.917\\
\noalign{\smallskip}\hline
\end{tabular}}
 \end{center}
  \end{minipage}
  \begin{minipage}[c]{0.4\hsize}
          \begin{center}
        \captionsetup{labelformat=empty,labelsep=none}
           \subcaption{$\beta=2, m=10,n=3$}
            \label{table1-b}
{\begin{tabular}{@{}cccc@{}} \toprule
$\alpha$&${F_\mathrm{sim}^{-1}}(\alpha)$ &$F_{40}^{-1}(\alpha)$ \\ \toprule
0.01	 &~0.451 & 0.451\\
0.05 &~0.539 & 0.539\\
0.50	 &~0.726 & 0.726\\
0.90	 &~0.833 & 0.834\\ 
0.95	 &~0.858 & 0.859\\ 
\noalign{\smallskip}\hline
\end{tabular}}
        \end{center}
    \end{minipage}
  \end {tabular}
  \end{center}
\end{table}
Fig~\ref{fig1} illustrates the probability $\Pr (0.7< \ell_n/\ell_1 < 1)$ with dimension for $\beta=1$ and $n=2$.
The probability for $\ell_n/\ell_1$ is large with the higher dimension; this implies that all eigenvalues are close to being equal under the null hypothesis. 
From (\ref{moment}), the mean, variance, skewness, and kurtosis of $1-\ell_n/\ell_1$ are provided in Table~\ref{table2}.
The mean and variance become smaller, while skewness and kurtosis become larger, as $m$ increases.
   \begin {table}[H]
  \caption {Mean, variance, skewness and kurtosis of $\ell_n/\ell_1$ for $n=2$ and $\beta=1$.}
  \label {table2}
  \begin{center}
\begin{tabular}{@{}cccccc@{}} \toprule
$m$ & Mean & Variance &  Skewness & Kurtosis & \\ \toprule
    5	 & 0.667 & 0.0415 & -0.675 & 2.81 \\
  25	 & 0.382 & 0.0246 & 0.0224 & 2.41\\
  45   & 0.301 & 0.0174 & 0.164 & 2.50 \\
  65   & 0.260 & 0.0136 & 0.235  & 2.57 \\
  85   & 0.232 &  0.0112 & 0.280 & 2.62 \\
 105  & 0.212 & 0.00961 & 0.311 & 2.66 \\
 125  & 0.196 & 0.00841 & 0.334 & 2.68 \\
 145  & 0.184 & 0.00749 & 0.349 & 2.70 \\
 \noalign{\smallskip}\hline
    \end {tabular}
   \end{center}
    \end {table}

\section*{Appendix. Recurrence relation of $q^\beta[\kappa, \mu]$ for $m=2, 3, 4$.}
\begin{itemize}
\item $m=2$, $\kappa=(\kappa_1, \kappa_2)$, $\mu=(\mu_1, \mu_2), \nu_1=\mu_1-\mu_2$.
\end{itemize}
\begin{align*}
\frac{\beta}{2}\{d^\beta(\mu)-d^\beta(\kappa)\}~q^\beta[\kappa, \mu] = (\nu_1+2)(\nu_1+1)~q^\beta[\kappa, (\mu_1+1, \mu_2-1)]
\end{align*}
\begin{itemize}
\item $m=3$, $\kappa=(\kappa_1, \kappa_2, \kappa_3)$, $\mu=(\mu_1, \mu_2, \mu_3), \nu_1=\mu_1-\mu_2, \nu_2=\mu_2-\mu_3$.
\end{itemize}
\begin{align*}
\frac{\beta}{2}\{d^\beta(\mu)-d^\beta(\kappa)\}~q^\beta[\kappa, \mu] =& (\nu_1+2)(\nu_1+1)~q^\beta[\kappa, (\mu_1+1, \mu_2-1, \mu_3)]\\
&+(\nu_2+2)(\nu_2+1)~q^\beta[\kappa, (\mu_1, \mu_2+1, \mu_3-1)]\\
&+3(\nu_1+1)(\nu_2+1)~q^\beta[\kappa, (\mu_1+1, \mu_2, \mu_3-1)]
\end{align*}\begin{itemize}
\item $m=4$, $\kappa=(\kappa_1, \kappa_2, \kappa_3, \kappa_4)$, $\mu=(\mu_1, \mu_2, \mu_3, \mu_4), \nu_1=\mu_1-\mu_2, \nu_2=\mu_2-\mu_3, \nu_3=\mu_3-\mu_4$.
\end{itemize}
\begin{align*}
\frac{\beta}{2}\{d^\beta(\mu)-d^\beta(\kappa)\}~q^\beta[\kappa, \mu] =& (\nu_1+2)(\nu_1+1)~q^\beta[\kappa, (\mu_1+1, \mu_2-1, \mu_3, \mu_4)]\\
&+(\nu_2+2)(\nu_2+1)~q^\beta[\kappa, (\mu_1, \mu_2+1, \mu_3-1, \mu_4)]\\
&+(\nu_3+2)(\nu_3+1)~q^\beta[\kappa, (\mu_1, \mu_2, \mu_3+1, \mu_4-1)]\\
&+2(\nu_2+2)(\nu_2+1)~q^\beta[\kappa, (\mu_1+1, \mu_2+1, \mu_3-1, \mu_4-1)]\\
&+3(\nu_1+1)(\nu_2+1)~q^\beta[\kappa, (\mu_1+1, \mu_2, \mu_3-1, \mu_4)]\\
&+3(\nu_2+1)(\nu_3+1)~q^\beta[\kappa, (\mu_1, \mu_2+1, \mu_3, \mu_4-1)]\\
&+4(\nu_1+1)(\nu_3+1)~q^\beta[\kappa, (\mu_1+1, \mu_2, \mu_3, \mu_4-1)]
\end{align*}




\begin{thebibliography}{00}


\bibitem{C1963}
A. G. Constantine, 
Some non-central distribution problems in multivariate analysis, 
{\it Ann. Math. Stat}. {\bf 34} (1963) 1270--1285. 


\bibitem{C1966}
A. G. Constantine, 
The Distribution of Hotelling's Generalised $T^2_0$, 
{\it Ann. Math. Stat}. {\bf 37} (1966) 215--225. 

\bibitem{Garcia2014}
J. A. D\'{i}az-Garc\'{i}a, 
Integral Properties of Zonal Spherical Functions, Hypergeometric Functions and Invariant, 
{\it J. Iran. Statist. Soc}. {\bf 13} (2014) 83--124. 

\bibitem{Garcia2011}
J. A. D\'{i}az-Garc\'{i}a and R. Guti\'errez-J\'aimez, 
On Wishart distribution: Some extensions, 
{\it Linear Algebra and Its Applications}. {\bf 435} (2011) 1296--1310. 

\bibitem{Dumitriu2007}
I. Dumitriu, A. Edelman and G. Shuman, 
MOPS: Multivariate orthogonal polynomials (symbolically), 
{\it Journal of Symbolic Computation}. {\bf 42} (2007) 587--620. 

\bibitem{Hashiguchi2000}
H. Hashiguchi, S. Nakagawa and N. Niki, 
Simplification of the Laplace-Beltrami operator, 
{\it Math. Comput. Simulation}. {\bf 51} (2000) 489--496.

\bibitem{Hayakawa1967}
T. Hayakawa, 
On the distribution of the maximum latent root of a positive definite symmetric random matrix, 
{\it Ann. Math. Stat}. {\bf 19} (1967) 1--17.

\bibitem{Jack1970}  
H. Jack,
Calculation and properties of zonal polynomials, 
{\it Proc. Roy. Soc. Edinburgh Sect. A}. {\bf 69} (1970) 1--18. 

\bibitem{James1964}  
A. T. James, 
Distributions of matrix variates and latent roots derived from normal samples, 
{\it Ann. Math. Stat}. {\bf 35} (1964) 475--501. 

\bibitem{Jiu2020}  
L. Jiu and C. Koutschan,
Calculation and properties of zonal polynomials, 
{\it Mathematics in Computer Science}. {\bf 14} (2020) 623--640. 

\bibitem{Kan2019}
R. Kan and P. Koev,
Densities of the extreme eigenvalues of Beta-MANOVA matrices, 
{\it Random Matrices: Theory Appl. } {\bf 8} (2019) 1950002. 

\bibitem{Khatri1968}
C. G. Khatri and K. C. S. Pillai,
 On the non-central distributions of two test criteria in multivariate analysis of variance, 
{\it Ann. Math. Stat}. {\bf 39} (1968) 215--226. 

\bibitem{Pillai1969}
 K. C. S. Pillai and T. Sugiyama,
Non-central distributions of the largest latent roots of three matrices in multivariate analysis, 
{\it Ann. Math. Stat}. {\bf 21} (1969) 321--327. 

 \bibitem{Shimizu2021a}
K. Shimizu and H. Hashiguchi, 
Heterogeneous hypergeometric functions with two matrix arguments and the exact distribution of the largest eigenvalue of a singular beta-Wishart matrix,  {\it J. Multivariate Anal}. {\bf 183} (2021) 104714.

\bibitem{Shimizu2021b}
K. Shimizu and H. Hashiguchi, 
Expressing the largest eigenvalue of a singular beta $F$-matrix with heterogeneous hypergeometric functions, 
{\it Random Matrices: Theory Appl}. (2021) 2250005. 

\bibitem{Shinozaki2021}
A. Shinozaki, S. Shimizu and H. Hashiguchi, 
Generalized heterogeneous hypergeometric functions and the distribution of the largest eigenvalue of an elliptical Wishart matrix, 
2021, arXiv:2104.12552.

\bibitem{Sugiyama1967}
T. Sugiyama, 
On the distribution of the largest latent root of the covariance matrix, 
{\it Ann. Math. Stat}. {\bf 38} (1967) 1148--1151. 

\bibitem{Sugiyama1970}
T. Sugiyama, 
Joint distribution of the extreme roots of a covariance matrix, 
{\it Ann. Math. Stat}. {\bf 41} (1970) 655--657. 

\bibitem{Takemura1984}
A. Takemura, 
Zonal Polynomials, 
Institute of Mathematical Statistics Lecture Notes-Monograph Series, 
Institute of Mathematical Statistics, Hayward (1984).

\bibitem{Ratnarajah2005} 
T. Ratnarajah and R. Vaillancourt, 
Complex singular Wishart matrices and applications, 
{\it Comput. and Math. Appl}. {\bf 50} (2005) 399--411.

\bibitem {Ratnarajah2004} 
T. Ratnarajah, R. Vaillancourt and M. Alvo, 
Eigenvalues and condition numbers of complex random matrices, 
{\it SIAM J. Matrix Anal}. {\bf 26} (2004) 441--456.

\bibitem{U1994}
H. Uhlig, 
On Singular Wishart and singular multivariate beta distributions, 
{\it Ann. Stat}. {\bf 22} (1994) 395--405. 

\end{thebibliography}


\end{document}